\documentclass{article}
\usepackage{amsfonts,amssymb,amscd, amsmath}
\usepackage{latexsym, graphicx}
\usepackage[affil-it]{authblk}
\usepackage{algorithmic}
\usepackage{algorithm}
\usepackage{tikz}
\usetikzlibrary{positioning}
\tikzset{main node/.style={circle,fill=blue!20,draw,minimum size=0.7cm,inner sep=0pt},
}

\usepackage{enumerate}

\newtheorem{theorem} {Theorem}[section]
\newtheorem{lemma}[theorem]{Lemma}

\newtheorem{conjecture}[theorem]{Conjecture}

\usepackage{xcolor}

\begin{document}

\title{On the structure of perfectly divisible graphs}
\author{Ch\'inh T. Ho\`ang }
\date{}
\maketitle
\begin{center}

Department of Computer Science and Physics\\ 
Wilfrid Laurier University\\ 
Waterloo, Ontario,  Canada, N2L 3C5\\
\texttt{choang@wlu.ca}

\end{center}

\begin{abstract}
A graph $G$ is perfectly
divisible if every induced subgraph $H$ of $G$ contains a set $X$
of vertices such that $X$ meets all largest cliques of $H$, and
$X$ induces a perfect graph. The chromatic number of a perfectly
divisible graph $G$ is bounded by $\omega^2$ where $\omega$
denotes the number of vertices in a largest clique of $G$. A graph $G$ is minimally non-perfectly divisible
if $G$ is not perfectly divisible but each of its proper induced subgraph is. 
A set $C$ of vertices of $G$ is a clique cutset if $C$ induces a clique in $G$, and $G-C$ is disconnected. 
We prove that a $P_5$-free minimally non-perfectly divisible graph cannot contain a clique cutset. This result allows us to re-establish several theorems on  the perfect divisibility of some classes of $P_5$-free graphs.
We will show that recognizing perfectly divisible graphs is NP-hard.  
\end{abstract}

\section{Introduction}\label{sec:introduction}
A {\em hole} is a chordless cycle with at least four vertices.  An
{\it antihole} is the complement of a hole. A hole is {\it odd} if
it has an odd number of vertices. A graph is odd-hole-free if it
does not contain, as an induced subgraph, an odd hole.
Odd-hole-free graphs have been studied in connection with perfect graphs.
%
%
For a graph $G$, $\chi(G)$ denotes the chromatic number of
$G$ and   $\omega(G)$ denotes the number of vertices in a largest clique of $G$.
A graph $G$ is {\it perfect} (\cite{Ber1961}) if for every
induced subgraph $H$ of $G$, we have $\chi(G) = \omega(G)$.
For more information on perfect graphs, the reader is referred to 
the books of Golumbic (\cite{Gol1980}), Berge and  Chv\'atal (\cite{BerChv1984}), 
Ramirez Alfonsin and Reed (\cite{RamRee2001}), and the survey paper of Ho\`ang and Sritharan (\cite{HoaSri2015}). 
For a set $S \subseteq V(G)$,
$G[S]$ denotes the subgraph of $G$ induced by the vertices of $S$.
Let the vertices of $V(G)$ be partitioned into two sets $A,B$. We say that $(A,B)$ is a {\it good} partition of $G$ if $G[A]$ is perfect, and $\omega(G[B]) < \omega(G)$. 
A graph $G$ is {\it perfectly divisible} (\cite{Hoa2018})  if every induced subgraph
$H$ of $G$ with at least one edge contains a good partition. The
chromatic number of a perfectly divisible graph $G$ is bounded by
$\omega^2(G)$. Ho\`ang  (\cite{Hoa2018}) conjectured that odd hole-free graphs are perfectly
divisible. 

\begin{figure}	
\begin{center}
\begin{tikzpicture}\label{fig:counter}
		\node[main node] (A) {$A$};
		\node[main node] (B) [right = 1cm  of A]  {$B$};
		\node[main node] (C) [right = 1cm  of B] {$C$};
		\node[main node] (D) [right = 1cm  of C] {$D$};
		\node[main node] (E) [below right = 2cm  of A] {$E$};
		\node[main node] (F) [right = 1cm  of E] {$F$};
		\node[main node] (G) [below  left = 2cm  of E] {$G$};
		\node[main node] (H) [right = 1cm  of G] {$H$};
		\node[main node] (I) [right = 1cm  of H] {$I$};
		\node[main node] (J) [right = 1cm  of I] {$J$};
		
		\path[draw,thick]
		(A) edge node {} (B)
		(B) edge node {} (C)
		(C) edge node {} (D)
		(E) edge node {} (A) 
		(E) edge node {} (F)
		(F) edge node {} (A)
		(F) edge node {} (D)
		(E) edge node {} (G)
		(E) edge node {} (J)
		(F) edge node {} (J)
		(G) edge node {} (H)
		(H) edge node {} (I)
		(I) edge node {} (J)
		;

\end{tikzpicture}
		\caption{A counter-example to the proof of Lemma 2.3 in \cite{DonXu2022b}. $A_1 =\{E,A,B,C,D\}, B_1 = \{F\}, A_2 = \{ F, J, I, H, G\}, B_2 = \{E\}$. The graph $G[A_1 \cup A_2]$ is not perfect. }
\end{center}
\end{figure}
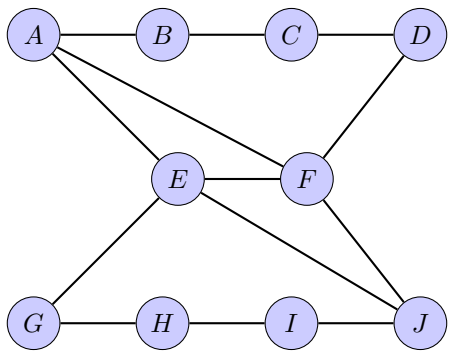

 A graph $G$ is {\it minimal imperfect} if $G$ is not perfect but each of its induced sugraphs is.
A graph $G$ is {\it  minimally non-perfectly divisible} (MNPD for short)
if $G$ is not perfectly divisible but each of its proper induced subgraph is.
A set $C$ of vertices of $G$ is a clique cutset if $C$ induces a clique in $G$, and $G-C$ is disconnected. The following result is folklore (for a proof, see \cite{Gol1980}).  
\begin{lemma}[Folklore]\label{lem:clique-cutset-perfect}
	No minimal imperfect graphs contain a clique cutset. 
\end{lemma}
It is natural to consider the analogy of Lemma \ref{lem:clique-cutset-perfect} for MNPD graphs. But this problem is still open. We will pose it as a conjecture.
\begin{conjecture}\label{conj:clique-cutset}
	No MNPD graphs contain a clique cutset.
\end{conjecture}
We note that Dong, Xu and Yu \cite{DonXu2022b} gave an incorrect proof of Conjecture \ref{conj:clique-cutset}. (In the proof of Lemma 2.3 in \cite{DonXu2022b}, the statement that $G[A_1 \cup A_2]$ is perfect is false. A counter-example is shown in Figure ~\ref{fig:counter}.) 

Let $C_k$ denote the hole on $k$ vertices. Let $P_k$ denote the chordless path on $k$ vertices. Let $G$ and $H$ be two graphs. We say $G$ is $H$-free if $G$ does not contain an induced subgraph isomorphic to $H$. Let $L$ be a set of graphs. A graph $G$ is $L$-free if $G$ does not contain
an induced subgraph isomorphic to a graph in $L$.

Dong, Xu and Yu proved the following ($K_{2,3}$ is the complete bipartite graph with one part having two vertices, and the other part having three vertices)
\begin{theorem}[\cite{DonXu2022b}]\label{thm:dong} 
	$(P_5, C_5, K_{2,3})$-free graphs are perfectly divisible.
\end{theorem}

The proof of Theorem~\ref{thm:dong} relies on the truth of Conjecture \ref{conj:clique-cutset} that as far as we know, is still unsolved. We will give a proof of the following
\begin{theorem}\label{thm:cutset}
	A $P_5$-free MNPD graph cannot contain a clique cutset.
\end{theorem}
The truth of Theorem~\ref{thm:cutset} means that Theorem~\ref{thm:dong} still holds.
We will also prove the following ($4K_1$ denotes the stable set on four vertices)
\begin{theorem}\label{thm:4K1}
	A $4K_1$-free MNPD graph cannot contain a clique cutset.
\end{theorem}
In Section \ref{sec:results}, we will establish a number of properties on MNPD graphs. In particular, we will prove Theorems \ref{thm:cutset} and \ref{thm:4K1}. In Section \ref{sec:recognition}, we will show that to recognize a perfectly divisible graph is a NP-hard problem. In Section \ref{sec:conclusions}, we will discuss open problems related to our work. 

\section{Properties of MNPD graphs}\label{sec:results}

Let $G$ be a graph. Recall that $\chi(G)$ denotes the chromatic number of
$G$, and  $\omega(G)$ denotes the number of vertices in a largest clique of $G$. The number of vertices in a largest stable set of $G$ is denoted by  
$\alpha(G)$. 
%
%
For a
graph $G$ and a set $X \subseteq V(G)$, $G[X]$ denotes the
subgraph of $G$ induced by $X$.


A vertex is {\it simplicial} if its neighbourhood is a clique. We will prove the following two results (for MNPD graphs) that are analogous to well known properties of minimal imperfect graphs. 
\begin{lemma}\label{lem:largest-clique}
	In a MNPD graph $G$, every vertex belong to a largest clique of $G$.
\end{lemma}
{\it Proof}. Let $G$ be a MNPD graph. Suppose $G$ contains a vertex $x$ that does not belong to a largest clique of $G$. Since the graph $G-x$ is perfectly divisible, its vertex set admits a good partition $(A,B)$. Let $B' = B \cup \{x\}$. Since $x$ belongs to no clique of size $\omega(G)$, we have $\omega(G[B']) < \omega(G)$. Now, $(A,B')$ is a good partition of $G$, a contradiction to the minimality of $G$. \hfill $\Box$

\begin{lemma}\label{lem:simplicial}
	No MNPD graphs contain a simplicial vertex.
\end{lemma}
\noindent
{\it Proof}. Let $G$ be a MNPD graph. Suppose $G$ contains a simplicial vertex $x$. Since the graph $G-x$ is perfectly divisible, its vertex set admits a good partition $(A,B)$.  Let $A' = A \cup \{x\}$. Since $x$ is simplicial, $G[A']$ is perfect, and so $(A',B)$ is a good partition of $G$, a contradiction to the minimality of $G$. \hfill $\Box$

We recall the Strong Perfect Graph Theorem of Chudnovsky, Robertson, Seymour, and Thomas (\cite{ChuRob2006}) that a graph is perfect if and only if it is (odd hole, odd antihole)-free. 

\noindent 
{\it Proof of Theorem ~\ref{thm:cutset}}. Let $G$ be a $P_5$-free MNPD graph. Suppose $G$ contains a clique cutset $C$. If $C$ is not a minimal cutset, then it contains a proper subset that is a minimal cutset of $G$. (Here, as usual, ``Minimal'' is meant with respect to set inclusion, not size). So, we may suppose that $C$ is a minimal cutset. Then $V(G) - C$ can be  partitioned into two non-empty sets $V_1, V_2$ such that there are no edges between $V_1$ and $V_2$. Let $G_i = G[C \cup V_i]$ for $i = 1,2$. By the minimality of $G$, for $i=1,2$, $G_i$ is perfectly divisible. So, for $i = 1,2$, $G_i$ admits a good partition $(A_i, B_i)$. 
Let $A = A_1 \cup A_2$, $B = (B_1 \setminus A_2) \cup (B_2 \setminus A_1)$. Then $(A,B)$ is a partition of $V(G)$. We will prove that the partition $(A,B)$ is good. 

Consider a largest clique $K$ of $G$. Then $K$ lies completely in $G_1$ or in $G_2$. But then $K$ must be met by $A_1$ or $A_2$, that is, $K$ must be met by $A$. So $G[B]$ contains no clique of cardinality $\omega(G)$, that is, $\omega(G[B]) < \omega(G)$. 

If $G[A]$ is perfect, then we are done. So we may assume that $G[A]$ contains as induced subgraph a graph $H$ that is an odd hole or odd antihole. Note that $G[A]$ has a clique cutset $C_A = (A_1 \cap C) \cup (A_2 \cap C)$. It follows that $H$ cannot contain a vertex in $V_1$ and a vertex in $V_2$. We may assume that $H \cap V_1 = \emptyset$, that is, $H$ lies completely in $G_2$. Since $G[A_2]$ is perfect, $H$ must contain a vertex in $A_1$. It follows that $H$ contains a vertex $a \in A_1 \cap C$. By the minimality of $C$, vertex $a$ has a neighbour $a' \in V_1$. If $H$ is an odd hole, then there are vertices $h_1, h_2, h_3 \in H$ such that $\{h_1, h_2 ,h_3\} \cap C = \emptyset$, and edges $h_1h_2, h_2h_3, h_3a$. 
But now $\{h_1, h_2 ,h_3, a, a'\}$ induces a $P_5$, a contradiction. Now we may assume $H$ induces an odd antihole of length at least 7 (note that the odd hole of length 5 is self-complementary). Consider the complement $\overline{H}$ of $H$. Since $C$ is a stable set in the complement $\overline{G}$ of $G$, $C$ does not contain two consecutive vertices of $\overline{H}$. Thus  $\overline{H}$ has vertices $h_1, h_2, h_3 \not\in C$ such that $\{h_1, h_2 ,h_3\} \cap C = \emptyset$, and edges $h_1h_2, h_2c, c h_3 \in E(\overline{G})$ for some vertex $c \in C \cap V(H)$. By the minimality of $C$, in $G$ vertex $c$ has a neighbour $c' \in V_1$. But now, $\{h_2, h_3, h_1, c, c'\}$ induces a $P_5$ in $G$, a contradiction. \hfill $\Box$

To prove Theorem \ref{thm:4K1}, we will need the following result.
\begin{lemma}\label{lem:one-side-perfect}
	Let $G$ be a graph with a clique cutset $C$.  Consider a partition of  $V(G) - C$  into two sets $V_1, V_2$ such that there are no edges between $V_1$ and $V_2$. If $G[C \cup V_1]$ is perfect, then $G$ is not MNPD. 
\end{lemma}
\noindent {\it Proof of Lemma \ref{lem:one-side-perfect}}. Let $G$ be a graph with a clique cutset $C$ and with $V_1, V_2$ defined as in the Lemma. Suppose that $G$ is MNPD. Let $G_i = G[C \cup V_i]$ for $i = 1,2$.  Suppose that  $G_1 = G[C \cup V_1]$ is perfect. Since $G_2$ is perfectly divisible, the vertices of $G_2$ can be partitioned into two sets $A_2, B_2$ such that $G[A_2]$ is perfect and $\omega(G[B_2]) < \omega (G_2) \leq \omega (G)$. Let $A = A_2 \cup V_1$, and $B = B_2$. We clearly have $\omega(G[B]) =\omega(G[B_2] < \omega(G)$. We will show that $G[A]$ is perfect. 
Consider an induced subgraph $G'$ of $G[A]$. Then $G'$ lies completely in $G_1$, or in $G[A_2]$, or contains a clique cutset that lies in $A_2 \cap C$. 
It follows from Lemma \ref{lem:clique-cutset-perfect} that $G[A]$ is perfect.  
But now $(A,B)$ is a good partition of $G$, a contradiction. $\Box$

\noindent
{\it Proof of Theorem \ref{thm:4K1}}. Let $G$ be a $4K_1$-free MNPD graph. Suppose $G$ contains a clique cutset $C$. Then $V(G) - C$ can be  partitioned into two sets $V_1, V_2$ such that there are no edges between $V_1$ and $V_2$. Let $G_i = G[C \cup V_i]$ for $i = 1,2$. 
If $\alpha(G[V_1]) > 1$ and  $\alpha(G[V_2]) > 1$, then $G$ contains a $4K_1$, a contradiction. So we may assume that $V_1$ is a clique. 
Then $G_1 = G[C \cup V_1]$ is the complement of a bipartite graph, and therefore a perfect graph.
But this is a contradiction to Lemma \ref{lem:one-side-perfect}. $\Box$

\section{The complexity of recognizing perfectly divisible graphs}\label{sec:recognition}
In this section, we will show that recognizing perfectly divisible graphs is a NP-hard problem. A graph $G$ is {\it triangle-free} if it does not contain a clique on three vertices, that is, $\omega(G) \leq 2$. Maffray and Preismann proved the following
\begin{theorem}[\cite{MafPre1996}]\label{thm:3-colorable}
	It is NP-hard to decide if a triangle-free graph is 3-colorable.
\end{theorem}
\begin{theorem}\label{thm:recognizing-perfectly-divisible}
	It is NP-hard to decide if a  triangle-free graph is perfectly divisible.
\end{theorem}
\noindent {\it Proof.} Let $G$ be a triangle-free graph. We only need to prove that $G$ is perfectly divisible if and only if $G$ is 3-colorable.  Without loss of generality, we may assume $G$ has at least one edge, and at least three vertices. Suppose $G$ is 3-colorable. Consider a 3-coloring of $G$. Let $A$ be the set of vertices with the first two colors, and $B$ be the set of vertices with the third color. Then we have  $\chi(G[A]) = 2$ and $\omega(G[B]) = 1 < \omega(G)$. Since $\chi(G[A]) = 2$, $G[A]$ is a bipartite graph, and therefore a perfect graph. So, $(A,B)$ is a good partition of $G$. 

Suppose that $G$ is perfectly divisible. Consider a good partition $(A,B)$ of $G$. Since $\omega(G[B]) < \omega(G) = 2$, $B$ must be a stable set. Since $G[A]$ is perfect and triangle-free, it is a bipartite graph. Thus, $\chi(G[A]) = 2$ and therefore $\chi(G) \leq \chi(G[A]) + \chi(G[B]) \leq 3$. The proof is completed. \hfill $\Box$

We do not know if perfect divisibility is a NP-property (for a discussion on NP-complete problems, see the classic text of Garey and Johnson \cite{GarJoh2002}).

\section{Conclusions and open problems}\label{sec:conclusions}

In this section, we discuss the connections between our results with other results and conjectures in the literature. 
Ho\`ang (\cite{Hoa2018}) proposed the following conjecture. 
\begin{conjecture}[\cite{Hoa2018}]\label{conj:odd-hole}
Odd hole-free graphs are perfectly divisible.
\end{conjecture}
We would like to propose a few more conjectures on perfectly divisible graphs. In connection to Theorem \ref{thm:cutset}, it is perhaps natural to propose the following

\begin{conjecture}\label{conj:P5-free}
	$P_5$-free graphs are perfectly divisible.
\end{conjecture}
Chudnovsky and Sivaraman \cite{ChuSiv2020} proved that $(P_5,bull)$-free graphs are perfectly divisible. 

In connection with Theorem \ref{thm:4K1}, we propose the following
\begin{conjecture}\label{conj:4K1}
	Graphs $G$ with $\alpha(G) \leq3$  are perfectly divisible.
\end{conjecture}
Ho\`ang (\cite{Hoa2018}, Lemma 7.3) proved that graphs $G$ with $\alpha(G) \leq 2$  are perfectly divisible. We note the conjecture of Sivaraman (\cite{KarKau2022}) that fork-free graphs are perfectly divisible (A {\it fork} is a graph with vertices $a,b,c,d_1, d_2$ and edges $ab,bc,cd_1, cd_2$). 

A hole is {\it even} if its length is even. Even hole-free graphs have been much studied in the literature (see  \cite{Vus2010} for a survey). A vertex is {\it bisimplicial} if its neighbourhood can be covered by two cliques, that is, the complement of its neighbourhood is a bipartite graph. 
Chudnosky and Seymour (\cite{ChuSey2023}) proved that every even hole-free graph contains a bisimplicial vertex. We would like to propose two conjectures.
\begin{conjecture}
	Even hole-free graphs are perfectly divisible.
\end{conjecture}
\begin{conjecture}
	No MNPD graphs contain a bisimplicial vertex. 
\end{conjecture}
 
A graph $G$ with at least one edge  is {\it $k$-divisible}  if the
vertex-set of each of its induced subgraphs $H$ with at least one
edge can be partitioned into $k$ sets, none of which contains a
largest clique of $H$. It is easy to see that the chromatic number
of a $k$-divisible graph is at most $k^{\omega -1}$. It was conjectured
by Ho\`ang and McDiarmid ({\cite{HoaMcd1999}, \cite{HoaMcd2002}})
that every odd-hole-free graph is 2-divisible. 
\begin{conjecture}[\cite{HoaMcd1999}]\label{conj:2-divisible}
	A graph is 2-divisible if and only if it is odd hole-free.
\end{conjecture}

Ho\`ang and McDiarmid proved the
conjecture for claw-free graphs (\cite{HoaMcd2002}). Ho\`ang (\cite{Hoa2018}) proved the conjecture for banner-free graphs.  
Dong, Song and Xu (\cite{DonXu2022}) proved
the conjecture for dart-free graphs. Chudnovsky and Sivaraman proved the conjecture for $(P_5, C_5)$-free graphs. The notion of clique cutset may also play an interesting role in the study of 2-divisible graphs. We would like to propose a conjecture analogous to Conjecture \ref{conj:clique-cutset}.
\begin{conjecture}
	No minimally non-2-divisible graphs contain a clique cutset.
\end{conjecture}

We will end the paper with the conjecture below.
 
\begin{conjecture}
	Even hole-free graphs are 3-divisible.
\end{conjecture}


\noindent {\bf Statements and Declarations}

The  author acknowledges the support of the Natural Sciences and Engineering Research Council of Canada (NSERC), 
[funding reference number DDG-2024-00015]. The author has no relevant financial or non-financial interests to disclose. No data were created or analyzed in this study. Data sharing is not applicable to this article.

\end{document}